\documentclass[12pt,thmsa]{article}
\usepackage{amsfonts}
\usepackage{amssymb}
\usepackage{graphicx}
\usepackage{amssymb,amsmath}
\newtheorem{teo}{Theorem}[section]

\newtheorem{lema}[teo]{Lemma}

\newtheorem{obs2}[teo]{Remark}

\newtheorem{tea}{Theorem}[subsection]

\newtheorem{no2}[teo]{Note}

\newtheorem{no3}[tea]{Note}

\newcommand{\SL}{{\rm SL}}
\newcommand{\Symm}{{\rm Symm}}

\newcommand{\PGL}{{\rm PGL}}

\newcommand{\GL}{{\rm GL}}

\newcommand{\F}{{\mathbb{F}}}
\newcommand{\Q}{{\mathbb{Q}}}

\def\timehm{\count31=\time \count32=\count31 \divide\count31 by 60
\number\count31 \multiply\count31 by 60 \advance\count32 by
-\count31 :\ifnum\count32<10 0\fi \number\count32}

\newcommand{\qed}{\hfill\rule{2mm}{2mm}}

\def\ideal#1{<\kern-2pt #1\kern-2pt >}

\begin{document}
\title{{\bf  Automorphy of $m$-fold tensor products of $\GL(2)$}
}
\author{Luis Dieulefait \thanks{Research partially supported by
 MTM2012-33830 and by an ICREA Academia Research Prize}
\\
Universitat de Barcelona\\
Departament d'Algebra i Geometria\\
G.V. de les Corts Catalanes 585\\
08007 - Barcelona, Spain \\
e-mail: ldieulefait@ub.edu\\
Telephone: 0034934021607\\
Fax: 0034934021601\\
 }
\date{\empty}

\maketitle

\vskip -20mm

\begin{abstract} We prove that for any $m > 1$ given any $m$-tuple of Hecke eigenforms $f_i$ of level $1$
 whose weights satisfy the usual regularity condition there is a
self-dual cuspidal
 automorphic form $\pi$
 of $\GL_{2^m}(\Q)$ corresponding to their tensor product, i.e., such that the system
 of Galois representations attached to
 $\pi$ agrees with the tensor product of the ones attached to the
 cuspforms $f_i$. \\

\end{abstract}

\newpage

\begin{quote}
\small \raggedleft \it -Eh bien, Axel, me dit mon oncle, cela va, et
le plus difficile
est fait. \\
-Comment, le plus difficile?  m'\'ecriai-je. \\
-Sans doute, nous n'avons plus qu'\`a descendre! \\
-Si vous le prenez ainsi, vous avez raison; mais enfin, apr\`es
avoir descendu, il faudra remonter, j'imagine? \\
-Oh!  cela ne m'inqui\`ete gu\`ere! \\
$ $ \\
  \rm ``Voyage au centre de la terre",   Jules Verne \\
\end{quote}

\section{Introduction}

This paper should be considered as a sequel to [10] since not only
we follow the strategy developed in [9] and [10] of propagating
automorphy through suitable ``safe chains" but also the skeleton of
this paper is the safe chain of $11$ steps carefully constructed in
[10]. Since we are going to reuse that chain (with just one minor
modification), the reader is warmly recommended to read its
construction in [10], section 3, before
reading this paper. \\

Our goal is to prove automorphy of $m$-fold tensor products of
Galois representations attached to level $1$ cuspforms, assuming
that such tensor product representations are regular (i.e., that
they have $2^m$ different Hodge-Tate numbers). We are going to use
the $11$ steps long chain from our previous paper (plus $4$ extra
moves at the end) in the first component (or factor) of the $m$-fold
tensor product in order to safely connect it with a CM form of prime
power level. At this point, automorphy of the tensor product will be
reduced to a similar automorphy statement for the $m-1$-fold tensor
product involving all but the first component. Thus, the result will
be proved by induction on $m$. This induction can of course start at
the trivial case $m=1$, or it can start at $m=2$ since for the case
of two newforms automorphy of the tensor product was proved by D.
Ramakrishnan (see [19]).
\\

There are, however, several technical problems that do not allow us
to do this in such a simple and straightforward way. Assuming that
in the process of moving the first component through the safe chain
 built in our previous
paper (while leaving the other components unchanged) the residual
images in the $m$-fold tensor product representations involved in
each congruence are sufficiently large for some Automorphy Lifting
Theorem (A.L.T.) to apply; which is something relatively easy to
guarantee by introducing Good-Dihedral primes in all components,
using information proved in [10] about the residual images in the
first component and applying some results about adequacy of tensor
products taken from [16]; there still remain two main
issues tricky enough to keep us busy for a while:\\
- the $m$-fold tensor products appearing in each congruence will not
necessarily be regular \\
- even if they were regular, the local conditions required to apply
some available A.L.T. would fail to hold oftentimes, for example,
the chain built in [10] requires working with small primes $p$, such
as $7$ or $11$, and for large $m$ and small $p$, assuming that the
$m$ representations in the tensor product are crystalline at $p$, it
is impossible to have $2^m$ different Hodge-Tate weights and to be
in a Fontaine-Laffaille situation at $p$, thus a priori it
is by no means clear that any A.L.T. can be applied. \\

In order to make the argument work in spite of the two serious
difficulties just described, we are going to develop two different
kind of tools: first, we are going to prove two variants of the
A.L.T. in [2] and in [11]. These will be A.L.T. for the special case
in which the two Galois representations involved are tensor
decomposable, one of them deals with the ``mixed case" of tensor
products where one factor is ordinary (not necessarily potentially
crystalline) and the others are potentially diagonalizable, and the
other is just an straightforward application of Harris' trick to the
case where one is just moving the weight of one component (assuming
potential diagonalizability in this component), leaving the others
unchanged (this kind of situation will appear very often in the
proof). Both theorems will be proved by following
the arguments in section 4 of [2]. \\
Secondly, as the reader may have guessed at this point, we are going
to be forced to make some non-trivial changes in the other
components (some ``genetic manipulations"), from the second to the last, of the tensor product. \\
In an analogy with Rubik's cube, you can think of this as if our
goal was to complete one face of the cube, yet in order to do so you
are forced to make some moves involving exclusively  the other
faces. \\
The congruences that we are going to apply on these other components
are thus not directed to link them to any sort of base case, they
are just meant to make all local conditions (plus regularity) to
hold in the $m$-fold tensor products while moving the first
component through the safe chain built in our previous paper. These
moves on the other components will thus be interspersed with the
moves on the first component, whenever required. Needless to say, in
any of these auxiliary congruences we must also make sure that some
A.L.T. applies  allowing automorphy to propagate. \\
The key tool that we will use in this process, in order to
manipulate the weights and behavior at $p$ of these other
components, is a result from [1] on existence of potentially
diagonalizable crystalline lifts of
arbitrary large weight for mod $p$ modular Galois representations.\\

This is a good point to mention that we will introduce two notions
on the weights of the modular forms appearing in the tensor
products, those of ``spread" and ``C-very spread" (where $C$ is a
positive integer). The latter being stronger than the former, both
notions imply regularity of the tensor product. We will use the tool
taken from [1] also to, at the very beginning, manipulate the
weights in order to reduce them to a situation where they are
``C-very spread" for a suitable value of $C$. This will be key to
guarantee regularity during the rest of the proof, in fact after
this step the weights will always be ``spread".\\

To summarize, our main task is to prove two useful variants of
A.L.T., and to explain the tricks in the just described
manipulations that are needed to make things work. There are other
issues, such as largeness and adequacy of residual images;
invariance of automorphy under twisting and Galois conjugation;
making sure that induction on $m$ work by, after having reduced the
problem to an $m-1$-fold tensor product, rewinding the process to
observe that such a tensor product is safely connected to the
initial one (with the first component dropped) whose components come
from $m-1$ level $1$ modular forms; etc. These and other issues will
be addressed on the go, we believe that none of them involves any
serious difficulty, they just need to be taken
care of. \\

Let us record the main result of this paper:\\

\begin{teo}
\label{teo:main}
Let $f_1, f_2, ...., f_m$ be an $m$-tuple of level $1$ cuspforms,
whose weights we will denote $k_1, ...., k_m$. We make the standard
regularity assumption that the $m$ integers $u_i = k_i -1$ not only
are all different but also the $2^m$ partial sums that one can form
with them are all different (we set $0$ as the sum of the elements
in the empty set).\\
Under this assumption, the $m$-fold tensor product of the
$\ell$-adic representations attached to the $f_i$ is automorphic and
the $m$-fold product $L$-function for the $f_i$ has analytic
continuation to the whole complex plane.
\end{teo}

Remarks: we are not going to insist in this paper on ``corollaries"
of this result (combined with the base change results from [9] and
[10]), we leave it as an exercise to show that for any totally real
number field $F$ unramified at $3$ one can lift this $m$-fold tensor
product to an automorphic representation on $\GL_{2^m}(F)$.
Moreover, changing the CM cuspform to be used as base case, one may
relax the assumption on the
field $F$.\\
Also, as happened in previous papers, conditional to some slight
strengthening of the available A.L.T., the theorem above (and its
base changed version) should also hold for newforms of higher level,
at least if the levels are odd (one can try to use weight reduction
to prove this unconditionally, but as far as we know not being
allowed to work in characteristic $2$ we can reduce the weight up to
$k= 4$ while it would be required to reduce it up to $k=2$). \\
Finally, for odd $p$, one can easily deduce that irreducible $m$-fold tensor
products of $2$-dimensional odd mod $p$ representations ramifying
only at $p$ are automorphic (in the sense that they admit an
automorphic $p$-adic lift) by combining the above
theorem with Serre's conjecture. Incidentally,  observe that
in this mod $p$ situation, the number of components $m$ must be
bounded as a function of $p$ in order to have irreducibility for the
tensor product (because there is only a finite number of level $1$
mod $p$
modular Galois representations). \\

Acknowledgments: I want to thank Toby Gee, Robert Guralnick and Richard Taylor for useful comments. In particular, I want to thank Robert Guralnick for providing the proof of adequacy for tensor products that we include in Section 3.

\section{In which two variants of the A.L.T. in [2] are proved}

Let us begin with the most difficult one, which is an A.L.T. for a
``mixed" case of an ordinary representation tensored with the tensor
product of some potentially diagonalizable ones. The proof does not
have any original idea in it, we are just going to follow closely
the proofs of Theorems 4.1.1 and 4.2.1 in [2] but with the
particular feature that all the important things, i.e., the
construction of a congruence ``against" an ordinary lift of large
weights and the application of Harris' trick to make weights (and
irreducible components at $\ell$) match, will be applied just to the
potentially diagonalizable components, and by just plugging in the
ordinary component in the two sides of the congruence we get what we
want. Since we need our result for ``small $\ell$" we have to work
as in [11] in a way that avoids the condition $\ell \geq 2 (n+1)$,
where $n$ is the dimension of the tensor product. Since it is good
enough for our applications, we will work for simplicity under the
condition $\ell \geq 2 (n_i + 1)$ for all $i$, where the $n_i$
denote the dimensions of the potentially diagonalizable components,
but see the remark after the proof of the Theorem on how to relax
this condition. The theorem that we can prove, in a fairly general
version good enough for our purposes, is the following (as usual,
A.L.T. are stated over imaginary CM fields, but a standard quadratic
base change argument provides a version of them over totally real
fields). \\

Given any $\ell$-adic  representation $x$ we denote by $\bar{x}$ the
semi-simplification of the corresponding residual representation.

\begin{teo}
\label{teo:mixed} Let $F$ be an imaginary CM field with maximal real
subfield $F^+$. \\
Let $\ell$ be an odd prime, and let \\
$$ (R, \mu_0), (S_i , \mu_i), i = 1, ...., t $$
be $t+1$ regular algebraic, irreducible, $n_i$-dimensional (with
$n_0 = \dim(R)$) polarized $\ell$-adic representations of $G_F$,
with:
$$ R \otimes  \bigotimes_i S_i $$
regular and irreducible. \\
Suppose that we are also given another $t+1$-uple of $\ell$-adic
representations of $G_F$ with the same dimensions $n_i$, all them
regular algebraic, irreducible and polarized, let us call them
$$ (R', \mu'_0), (S'_i , \mu'_i), i = 1, ...., t $$
On the level of residual representations, we suppose that $(R,
\mu_0)$ and $(R', \mu'_0)$ are congruent mod $\ell$, and similarly,
for every $i$, we assume that $(S_i , \mu_i)$ and $(S'_i , \mu'_i)$
are congruent mod $\ell$.  \\
We also assume that the tensor product
$$ R' \otimes  \bigotimes_i S'_i $$
is automorphic, regular algebraic, irreducible and polarized. We
make the following technical assumptions:\\
1) For every place $v$ of $F$ dividing $\ell$, $R|_{G_{F_v}}$ and
$R'|_{G_{F_v}}$ are ordinary and $S_i|_{G_{F_v}}$ and
$S'_i|_{G_{F_v}}$ are potentially diagonalizable, for every $i$.\\
2) $\overline{R \otimes  \bigotimes_i S_i   } (G_{F(\zeta_\ell)}) $
is
adequate, and $\zeta_\ell \notin F$.\\
3) For every $i > 0$, $\ell \geq 2 (n_i + 1)$.\\
Then:
$$ R \otimes  \bigotimes_i S_i $$
is automorphic.
\end{teo}

Proof: the proof consists on applying the same procedure used in
Theorems 4.1.1 and 4.2.1 in [2], ignoring for the most part the
ordinary component $R$ (or $R'$), just plugging in $R$ ($R'$,
respectively) after Harris' trick having been applied to the
remaining components: this works fine because the ``connectedness"
relation from [2] (at any place $v$) behaves well under tensor
products, and a representation is clearly connected to itself.
Finally, for the use of Geraghty's theorem,  tensoring with an
ordinary representation of sufficiently spread weights preserves
ordinariness, as we will shall see, thus we will be comparing
ordinary representations at the key step of changing weights and
types, as required. We limit ourselves to mention the few
differences, but we are not going to reproduce here the application
of Harris' trick since it is to be
applied exactly as in Theorem 4.1.1 of loc. cit. \\
It is important to remark that condition (2) implies that for every
$i > 0$ $\bar{S}_i((G_{F(\zeta_\ell)}))$ is irreducible  (*). \\
Construction of the ordinary lifts: as in Theorem 4.2.1 of loc. cit,
this has to be applied twice, one of each side, so to ease notation
let us do it for the $S_i$ (a similar procedure applies to the
$S'_i$). As a general rule, when we go to a solvable extension, even
if the construction only involves the $S_i$, we are going to put the
extra condition that such an extension is linearly disjoint from the
fixed field of $\ker{\bar{R}}$, so that property (2) is preserved on
such an extension. With just this extra restriction in mind, here we
proceed as in loc. cit.: we take a solvable extension $F'$ of $F$
where locally at all places above $\ell$ all the $\bar{S}_i$ have
trivial image, and there by Theorem 3.2.1 of loc. cit. (which can be
applied thanks to condition (3) above and the irreducibility
condition (*)) we can construct a crystalline (at each place above
$\ell$) ordinary lift of $\bar{S}_i$ of arbitrary regular Hodge-Tate
weights, such that it is connected to $S_i$ locally at every place
$v$ not dividing $\ell$. If we call $W_i$ these ordinary lifts, in
order to have regularity later when applying Harris' trick we need
the following to hold:
$$ R \otimes \bigotimes_i  S_i \otimes \bigotimes_i  W_i  \qquad \quad (**)$$
is regular. Let us show how to pick the weights of each $W_i$ in
order that this condition is satisfied, using an iterated version of
what is done in loc. cit., Theorem 4.2.1. For the ordinary lift
$W_1$, we pick an integer $t_1$ such that it is bigger than $|h  -
h'  |$ for any pair $h, h'$ of different Hodge-Tate weights of $ R
\otimes \bigotimes_i S_i $, and we choose as weights: $\{ 0, t_1,
.... (n_1 -1) t_1  \}$.
\\
Having done this for all $W_i$ for $i < j \leq t$, for the ordinary
lift $W_j$ we pick an integer $t_j$ such that it is bigger than $|h
- h'  |$ for any pair $h, h'$ of different Hodge-Tate weights of $$
R \otimes \bigotimes_i S_i \otimes \bigotimes_{i < j}  W_i $$ and we
choose as weights: $\{ 0, t_j, .... (n_j -1) t_j  \}$ \\
This concludes the construction of the ordinary lift of the tensor
product of the $S_i$, which is the tensor product of the $W_i$.
Moreover, and this will be key later in order to apply Geraghty's
theorem, $ R \otimes  \bigotimes_i W_i $ is ordinary. \\
At the
referee's request, let us expand on this: tensoring two ordinary
representations, even when we know that the tensor product is
regular,
 is not enough to guarantee that the tensor product is ordinary. The condition that may fail
  is that the Hodge-Tate weights shall be increasing
  as we move through the graded pieces corresponding to the filtration by Galois invariant subspaces
  (see for example [2] for the precise definition of ordinariness). But an easy computation shows inductively that
  in the tensor product under consideration
   the Hodge-Tate weights of each $W_i$ have been chosen to be sufficiently spread so that not only regularity
    but also ordinariness
    of the tensor product holds.  In fact, the easier way to perform this computation is by considering the tensor products
     in the opposite sense: $W_1 \otimes R$, $W_2 \otimes W_1 \otimes R$, etc. The required condition on the increasing of the Hodge-Tate
      weights along the graded pieces follows then automatically from the inequalities imposed on the $t_i$ (and, of course,
      from the assumption that each representation in the tensor product is itself ordinary).
      For example, $W_1 \otimes R$ is ordinary because the difference of two consecutive weights of $W_1$ is $t_1$ and by
      construction $t_1$ is larger than the difference between the largest and the smallest weight of $R$.
       This proves that the tensor product is ordinary, by
       just observing that in an irreducible  tensor product permuting the order of the factors
       gives always the same representation, up to isomorphism (alternatively, one can prove first
       the A.L.T. with $R$ and $R'$ as last factors on both sides,
       and again because permuting the order of the factors in a tensor product gives isomorphic representations
       the result follows as stated since automorphy of a Galois representation is preserved by an isomorphism).  \\
       \newline
Here, as in Theorem 4.1.1 of loc. cit, Harris' trick is applied to
``force" $\bigotimes_i S_i$ and $\bigotimes_i W_i$ to connect also
locally at places dividing $\ell$ by tensoring with certain
well-constructed monomial representations, both with the same
residual image. Observe that this makes sense because we have
assumption (1) implying that $\bigotimes_i S_i$   is potentially
diagonalizable. At this point, we tensor by $R$ both sides,
something that preserves connectedness at every finite place.
Regularity is ensured in this process (equivalently, (**) is
regular) thanks to our choice of the weights of the ordinary lifts
$W_i$. As for the images: we stress once again that in all auxiliary
fields, including the field $M$ where the characters giving rise to
the monomial representations in Harris' trick are defined, we have
to incorporate the extra condition that they should be linearly
disjoint from the fixed field of the kernel of $\bar{R}$. This being
said, we have from assumption (2) that the residual images of the
tensor product of the $R$ and the $S_i$'s will be adequate (even
when restricted to the cyclotomic extension) and by applying Lemma
A.3.1 in [1] we see that tensoring with a well-suited monomial
representation (i.e., one satisfying the linear disjointness
condition just discussed) will preserve this adequacy of the
residual images. \\
After playing Harris' trick, we obtain two tensor product
representations and as in section 4 of loc. cit., we would like to
apply Thorne's
A.L.T. (Theorem 2.3.1 in loc. cit.) to conclude this step. \\
From what we have done we see that all conditions needed to apply
Thorne's Theorem are satisfied, except for the fact that we may be
dealing with representations that are not potentially crystalline,
since $R$ may not be so.\\
Luckily, as observed by Calegari in [6], Thorne's Theorem easily
extends to the case of potentially semistable representations. When
written, this observation in [6] required two technical conditions
(smoothness of the point in the local deformation ring at places
above $\ell$, and Shin-regularity of the weights), but since then
local-to-global compatibility at $\ell =p $ has been fully
established by Caraiani in [7], and it is easy to see that as a
consequence the variant of Calegari to Thorne's result holds in
general (the smoothness condition is always satisfied as follows
from the proof of Lemma 1.3.2 of [2], see [6], Remark 2.8; and the
Shin-regularity, only needed in [6] to use a geometric model of the
Galois representations attached to automorphic forms, which in turn
is just used to deduce finiteness of the number of possible types at
$p$ that show up during the Taylor-Wiles process, is not necessary
now since it follows clearly from local-to-global compatibility that
by fixing the level at $p$ the  types
 locally at $p$ of the $p$-adic Galois representations that occur
 are finite in number, a finite set of types independent
  of the auxiliary primes introduced to the level of the automorphic forms appearing in the
  Taylor-Wiles construction). \\
  This being said, we apply this variant of Thorne's Theorem and the
  outcome is that the automorphy of $ R \otimes  \bigotimes_i S_i $
  follows from that of $ R \otimes  \bigotimes_i W_i $. Here we are hiding the core idea of Harris' trick, namely, that
  automorphy behaves well, in both directions, under tensoring by a monomial representation (see [2] for details). \\
  The rest of the argument is verbatim as in Theorem 4.2.1 of
  [2]: all that remains is a second application of the above
  procedure, on the other side, and then an application of
  Geraghty's Theorem ensuring that automorphy propagates well among
  ordinary representations, which we can apply due to assumption
  (2), and the care we have had in choosing the solvable extensions.
   This concludes the proof.

\qed\\

Remark: Condition (3) can be replaced by the condition: \\
(3') for every $i$, there is a finite extension $F'_i$ of $F$ such
that after restriction to $G_{F'_i}$ the residual representation
$\bar{S}_i$ is potentially diagonalizably automorphic and $\bar{S}_i
(G_{F'_i(\zeta_\ell)})$ is adequate.\\
In fact, using Theorem A.4.1 in [1] instead of Theorem 3.2.1 in
[2] the proof follows with the same arguments.\\
(3') is strictly weaker than (3) when combined with the other
conditions in the theorem: assuming (3), since we also have
condition (2) implying that $\bar{S}_i((G_{F(\zeta_\ell)}))$ is
irreducible for every $i
>0$, we can apply Corollary 4.5.3 of [2] together with the main
result in [16]
to deduce that (3') also holds. \\

The second A.L.T. that we propose is a variant of Theorem 4.1.1 in
[2] with the addendum that
 tensor product by a common representation (or by two representations that are connected locally
  everywhere, one on each side) has been taken. The statement is the following:

\begin{teo}
\label{teo:enlamenor}
Let $F$ be an imaginary CM field with maximal real
subfield $F^+$. \\
Let $\ell$ be an odd prime, and let \\
$$ (R, \mu_0), (S, \mu_1)$$
be two regular algebraic, irreducible,  polarized $\ell$-adic representations of $G_F$,
with:
$ R \otimes  S$
regular and irreducible. \\
Suppose that we are also given another    pair of $\ell$-adic
representations of $G_F$
$$ (R', \mu'_0), (S', \mu'_1)$$
of the same dimension than $R$ and $S$, respectively, both of them
regular algebraic, irreducible and polarized, and such that $R' \otimes S'$ is regular and irreducible. \\
On the level of residual representations, we suppose that $(R,
\mu_0)$ and $(R', \mu'_0)$ are congruent mod $\ell$, and that $(S , \mu_1)$ and $(S' , \mu'_1)$
are congruent mod $\ell$.  \\
We also assume that the tensor product
$$ R' \otimes  S'$$
is automorphic. We
make the following technical assumptions:\\
1) For every place $v$ of $F$ dividing $\ell$, $R|_{G_{F_v}}$ and
$R'|_{G_{F_v}}$ are potentially diagonalizable, and for every $v$ not dividing $\ell$,
$R|_{G_{F_v}} \sim R'|_{G_{F_v}}$. \\
2) $S$ and $S'$, locally at any finite place $v$, are connected.\\
3) $\overline{R \otimes S} (G_{F(\zeta)}) $ is adequate.\\
4) $ R \otimes R' \otimes S$ is regular. \\
Then, $R \otimes S$ is automorphic.
\end{teo}

The proof is a straightforward adaptation of the proof of Theorem
4.1.1 in [2],
 and the modifications required are included in what we explained during the proof
 of the previous theorem. For the reader convenience, let us repeat the relevant
 arguments: we proceed as in Theorem 4.1.1 of loc. cit., applying Harris' trick to
  $R$ and $R'$. Through the construction of a suitable pair of monomial representations,
   over a suitable extension, and tensoring with them $R$ and, respectively, $R'$, we manage
    to obtain a congruence between representations that now connect locally at every place. This
     is possible thanks to assumption (1). If we now tensor one side by $S$ and the other by $S'$, due
      to assumption (2), we still obtain a congruence between two representations that connect locally at
       every place. At this point, we want to apply Thorne's Theorem (2.3.1 in loc. cit.), or maybe
        Calegari's variant to it (in case $S$ is not potentially crystalline at some $v$ dividing $\ell$)
         recalled during the proof of the previous theorem, to these two triple tensor products. So it
          remains to check that the conditions required to apply this A.L.T. hold. As in the course of the previous proof, one should take the extra care during the applications of Harris' trick to choose all solvable extensions and auxiliary number fields to be linearly disjoint from the fixed field of the kernel of $\bar{S}$. In particular, if the image of $\bar{S}$ is monomial, we make sure that the monomial representations constructed in Harris' trick are induced from  characters of a field that is linearly disjoint from the one from which  $\bar{S}$ is induced. This being said, condition (3) together with Lemma A.3.1 in [1] guarantee the required adequacy condition on residual images. Regularity is ensured by condition (4), which is the natural generalization of a similar condition in Theorem 4.1.1 of [2]. This implies that the A.L.T. of Thorne can be applied and we conclude as in Harris' trick that $R \otimes S$ is automorphic because automorphy behaves well under tensor products with monomial representations. \\
Remark: as in [11], we have applied results from [1] to be able to prove a theorem valid also for ``small" primes. \\

\section{In which some important features of the safe chain from [10] are recalled,
 and a little extra effort is done to end the chain in a CM form}

We begin by reproducing the description taken from section 2 of [10]
of the $11$ steps long safe chain. In all the congruences appearing
in the following $11$ steps residual images are ``large", as
 defined in the Introduction to loc. cit. As proved in [15] this implies that these images are adequate,
  even after restriction to any solvable extension. Recall also that all these congruences are in
   characteristic larger than $5$.\\
We start with a cuspform $f$ of level $1$ and weight $k \geq 12$: \\

1) Introduce Good-Dihedral prime (as in [18]): level raises to $q^2$. \\
2) Weight reduction via Galois conjugation: weight reduces to $2 < k \leq 14$ ($k$ even).\\
 3) {\it Ad hoc} tricks to
make the small weight congruent to $2$ mod $3$ (Sophie Germain
primes, Hida families...): end up with $k \equiv 2 \pmod{3}$ and $k
< 43$.
\\
4) Introduce MGD prime $43$ using the pivot primes $7$ and $11$ (as
in [9]): end up with newform of weight $2$ and
level $43^2 \cdot q^2$.\\
5) Remove the Good-Dihedral prime (in two moves): end up with a newform of weight $q+1$ and level $43^2$.\\
6)  Again weight reduction via Galois conjugation but this time
``highly improved", because we need to ensure large residual image
at each step using just the MGD prime.
Weight reduces to  $2 < k \leq 14$ ($k$ even).\\
7) {\it Ad hoc} tricks to make weight smaller than $17$ and divisible by $4$
 (Sophie Germain primes, Khare's weight
reduction...): end up with a newform of weight
 $16$ and level $43^2$. \\
8)  Introduce nebentypus at $17$ of order $8$: end up with a newform of
 weight $2$, level $43^2 \cdot 17$, with nebentypus. \\
9)  Remove the MGD prime modulo $11$ via an {\it ad hoc} Lemma to
ensure residual irreducibility: get congruence (maybe using
level-raising)
 with a newform Steinberg at $43$, of weight $2$ and level $43 \cdot 17$,
 with nebentypus of order $8$ at $17$. \\
10) Move from weight $2$ to weight $44$ by reducing modulo $43$: irreducibility checked by hand. \\
11)  As predicted by Generalized Maeda in weight $44$, level $17$, nebentypus of
order
$8$: check that this space has a unique orbit.\\

For what we will do later, it is easier to introduce a
simplification in Step 1: in [10], before performing the level
raising at $q$, we first make a move to reduce to a weight $2$
situation.
 Here we prefer to do it directly in the given weight: if $f$ is a level $1$ newform of weight $k$, we
 select a bound $B$ larger than $68$ and larger
  than $2k$ (these two inequalities appear also in Steps 1 and 2 of loc. cit., they are required at some steps), and by taking a prime $t$ larger than $B$ and congruent to $1$ mod $4$ such that the image of the mod $t$ Galois representation attached to $f$ is $\GL(2, \F_t)$ one can through level raising add a Good-Dihedral prime $q$ (Good-Dihedral with respect to the bound $B$), a supercuspidal prime in the level, satisfying $t \mid q+1$, $q \equiv 1 \pmod{8}$,
 and $q \equiv 1 \pmod{p}$ for every $p \leq B$. The reader should consult [12] where introduction of Good-Dihedral primes (called there ``tamely dihedral'') for newforms of weight $k>2$ is discussed in full detail. Ramification at $q$ of the level-raised Galois representation is given by a character of order $t$ of the quadratic unramified extension of $\Q_q$. The key point of a Good-Dihedral prime is that as long as characteristics $t$ and $q$ are avoided, the residual images will be absolutely irreducible, and as long as one works in characteristic $p \leq B$, or any characteristic $p$ such that $q$ is a square mod $p$ (which, due to quadratic reciprocity, is equivalent to require that $p$ is a square mod $q$) the residual images will be large. \\
In fact, if we rename now $f$ to $f_1$ and
 $k, B, t,q$ to $k_1, B_1 , t_1, q_1 $, where
 $f_1$ is the newform of level $1$ of smallest
 weight $k_1$ in a regular (in the sense used in Theorem \ref{teo:main})
  $m$-tuple of newforms $f_i$ of weights $k_i$ (we suppose for simplicity that
  the sequence $k_i$ is increasing), the very first move that we will do to prove
  Theorem \ref{teo:main} will be this level raising in the first component $f_1$
  modulo $t_1$. Thus, it is very convenient to impose the following extra conditions:
   when choosing $B_1$, take it also to be larger than all $k_i$, for all $i$, thus
   since $t_1 > B_1$ when reducing modulo $t_1$ the $m$ residual representations will
   all be Fontaine-Laffaille, thus potentially diagonalizable, at $t_1$. Also, take $t_1$
    so that the $m$ components have all large, therefore adequate, residual images (this is possible due to
    the large images Theorem of Ribet, see [20]). It is worth stressing that the
    definition of large (see [10], Introduction) implies
     that the projectivized residual images will be (non-abelian) almost simple groups. Finally, observe that since we are
    in a Fontaine-Laffaille situation and all weights are different, it is easy to see by
    looking at the action of the inertia group at $t_1$ that the $m$ residual representations
 are  not isomorphic to each other, not even up to twist, and this implies from the results in
 [16] that the tensor product of the $t_i$-adic Galois representations attached to the $f_i$
  has residually adequate image\footnote{let us elaborate on this:
 if we take the definition of adequacy from [1], Appendix A,
 we see trivially that conditions (1) and (2) are satisfied,
 condition (3) follows from the fact that the components satisfy this condition as proved in [16], Lemma
 2.(ii), and finally for condition (4) applying K\"unneth's formula
 one can see again that it is satisfied because it is satisfied by
 the components (just write the adjoint representation as the tensor product of the representation with its dual).
  We thank R. Guralnick for providing this proof}. The same is true if we restrict to  any solvable extension,
   since largeness is preserved by doing such thing. Therefore, we see that we have a
   congruence modulo $t_1$ where the good-dihedral prime $q_1$ is being introduced to
   the first component, and the other components are left unchanged, such that Theorem 4.2.1 in
    [2] can be applied to show that automorphy of one side is equivalent to automorphy of the other. \\

After this simplified version of Step 1, the 11 steps go exactly as in Section 3 of [10]. This is what will happen in the first component $f_1$, we have seen in the previous paragraph how Step 1 can be done in such a way that an A.L.T. applies to the $m$-fold tensor product of the $f_i$, the rest of this paper will be devoted to either just observe or make suitable manipulations in the other components to make this happen, that we can go along the 11 steps in the first component in such a way that some A.L.T. can be applied at any of the involved congruences between $m$-fold tensor products. \\

But before getting into this, it remains to complement the chain (we
want it to end in a CM point): as it is it ends, in Step 11 above,
 in a  space of newforms of ``small" level and weight with a single orbit
  of conjugated newforms, we will call this space the ``bottom space". We
  can try to argue as in the  proof of base change in [10] that we can
  use the exact same strategy to connect an unspecified  CM form to this bottom space, but
  there we were working exclusively in dimension $2$, where very powerful A.L.T.
  of Kisin made the level-reduction possible. To remedy this we will show by direct
  computations that a suitable CM form can be connected to the bottom space in just
  four moves, and the local conditions showing up in each of these four congruences
   will be similar to those already appearing in the 11 steps above, so it is a safe
    way to complement our chain. We will explain the four moves, just remember that
     they have to be applied in reverse order when linking this to our chain,
     starting from a newform in the bottom space, and ending in the CM form. \\
The CM form $g_1$ that we consider has weight $2$ and level $27$, it
corresponds to a CM elliptic curve. We check that to this newform,
the prime $43$ can be added as a Steinberg prime to the level modulo
$13$, in fact its eigenvalue satisfies $a_{43} = 8 \equiv \pm 44
\pmod{13}$ (i.e., the condition under which
 Ribet proved level-raising holds). Moreover, we do a few computations in the
 space of newforms of weight $2$ and level $43 \cdot 27$ and we find the
  newform it is congruent to:  a non-CM newform $g_2$ whose field of coefficients $\Q_{g_2}$ has degree $10$
  having a dihedral mod $13$ Galois representation
   (dihedral but no bad-dihedral, since it is induced from $\Q(\sqrt{-3})$). This newform is listed $13$-th and
   last in William Stein's table of newforms of level $1161$, weight $2$
   and trivial nebentypus, see
    http://modular.math.washington.edu/Tables/Eigenforms/eig1101-1200 .
   \\
   The $3$-adic Galois representations attached to $g_2$ are potentially
Barsotti-Tate
     and there is a prime dividing $3$ in $\Q_{g_2}$ whose
      inertial degree is $3$, and we compute the corresponding residual representation
       and see that its image contains $\SL_2(\F_{27})$ (in fact the projective image is
       $\PGL_2(\F_{27})$)
       and it is thus large and
       adequate (see [15]). We check that, for the chosen prime above $3$ and residual representation,   there are no congruences (not even
        up to twist) with newforms of level $43$ and weight $2$, and there is a
         congruence with one of level $43$ and weight $4$, which we call
         $g_3$. This newform has a field of coefficients of degree $6$, it is listed second and last in the table of newforms of the corresponding space,
         see
         http://modular.math.washington.edu/Tables/Eigenforms/eigk4-1-100. \\
          Since the mod $3$ congruence between $g_2$ and $g_3$ will be part of our chain,
           we need to know more about $g_3$ locally at $3$. It is crystalline, because $3$
            is not in the level of $g_3$, but it fails to be Fontaine-Laffaille. Fortunately,
             since we know (by computations) that there is no congruence with a newform of
             weight $2$ and level $43$, this implies that $g_2$ must be ordinary at $3$ (to
              see this, combine the results in [3] with the strong form of Serre's conjecture).
               Thus it is ordinary and crystalline, hence potentially diagonalizable. \\
Now we do the last two moves, which are standard: we reduce mod $17$ and take a weight $2$ lift,
 whose nebentypus is $\omega^2$, of order $8$, and then, as in Step 10 above, reduce mod $43$
 to end up in the bottom space. The last step is thus exactly as Step 10 and irreducibility and
  largeness of the residual image follow as in that step (see [10]). In the previous step,
  working modulo $17$, with weight $4$ and level $43$, we have a congruence between a
  Fontaine-Laffaille representation and a potentially Barsotti-Tate one, and the only
  problem is that we should check that the residual image is large. But this follows
  from the results in [4] (see ``Theorem for square-free level case", in the Introduction
   of that paper) since $17 \nmid 43^4 - 1$ and $17 > 4k - 3 = 13$. This completes the
    four moves required to go from the bottom space to a CM form, or viceversa. \\

{\bf Definition}: The safe chain that we have just recalled, and complemented, starting
on the system of Galois representations attached to the level $1$ newform $f_1$ and
finishing in the system attached to a $CM$ form of weight $2$ and level $27$ will be called the {\it skeleton chain}. \\

\section{In which spread and $C$-very spread weights are defined, and the proof is reduced to a $C$-very spread case}

Recall that, by assumption, in Theorem \ref{teo:main} we have $m$ newforms $f_i$ whose weights $k_i$ are such that the tensor-product Galois representation is regular. In such a case, we say that the $m$-tuple of weights $0< k_1 < k_2 < ...... < k_m$ is  {\it regular}. \\
We make two more restrictive definitions on an $m$-tuple of weights:\\

{\bf Definition} Given an $m$-tuple of weights $k_1, k_2,...., k_m$ we say that it is {\it spread} if the following inequalities hold: \\
$$ k_2  >  k_1 $$
$$ k_3  >  k_1  + k_2 $$
And, in general, for every $2 \leq i \leq m$:
$$k_i > k_1 + k_2 + ...... + k_{i-1} $$

{\bf Definition} Suppose that we have an $m$-tuple of weights $k_1 <
k_2 < ......< k_m $ such that the first weight $k_1$ is a variable
(i.e., it is not constant), and it is constrained to move in the
interval
$$ 2 \leq k_1  \leq  C   \quad \quad (*)$$
for certain constant $C >2 $. Then, we say that the $m$-tuple is {\it $C$-very spread} if, together with the condition (*) in the varying first weight, the other weights are constant and satisfy:
$$  k_2  >  2 C $$
$$ k_3  >   2 C + k_2$$
And, in general, for every $2 \leq i \leq m$:
$$k_i > 2 C + k_2 + ...... + k_{i-1} $$

Obviously, if $k_1$ is fixed and satisfies (*) for some $C$, being $C$-very spread implies being spread. Also, spread implies regular, as it is easy to check.\\
If we forget for one second about the factor $2$ in the $2 C$ appearing in
the definition of $C$-spread, the reason for this notion is clear: we will
be moving the first weight in a certain range, bounded above by certain $C$,
and we want the set of weights to be spread along the process, therefore regular.
 If we know the value of $C$ and if we assume that the weights $k_2, ...., k_m$ are
  as in the definition of $C$-spread, then this means that the weights will be spread
  even in the worst case where $k_1 = C$, thus we win: we have regularity through the whole process. \\
The reason for the $2$ multiplying $C$ is that sometimes we will be
moving $k_1$ (the weight of a newform showing up in the first
component of our $m$-fold tensor product) always taking values
bounded by some $C$, and thus we will have a congruence, in the
first component, between two newforms $f_1$ and $f'_1$ of different
weights $k_1$ and $k'_1$, both bounded by $C$. Then, as it will
happen many many times in our proof, suppose that we want to apply
Theorem \ref{teo:enlamenor} to this congruence of $m$-fold tensor
products, i.e., with the Galois representations attached to $f_1$
and $f'_1$ playing the role of $R$ and $R'$ and with $S = S'$
corresponding to the Galois representation  which is the $m-1$-fold
tensor product of the other components. If you assume that all
conditions except number (4) of that Theorem are satisfied (in
particular, this means that we must be in a situation where both
$f_1$ and $f'_1$ are potentially diagonalizable), still in order to apply it we have
to make sure that the ``extra-regularity" condition (4) is
satisfied, i.e., if we say (just to ease notation) that $k_1 < k'_1$
we want the $m+1$-tuple of weights:
$$(k_1 , k'_1 , k_2 , k_3 ,......, k_m) $$
to be regular. But $k_1 < k'_1 \leq C$ together with the other inequalities in the definition of $C$-very
 spread imply that this $m+1$-tuple is spread, thus regular, so this extra-regularity condition is satisfied. \\
This is perhaps a good moment to reveal what will the value of $C$ be in our proof: as
 some readers may have already guessed, we will take $C = q_1 +1$, where $q_1$ is the
 Good-Dihedral prime already introduced to the first component. Looking at the skeleton chain from section 3, it is easy to see
 that as the first component moves along this chain, its weight $k_1$ satisfies
 $k_1 \leq q_1 + 1$. From now on, we fix:  $C : = q_1  + 1$.\\

But the starting situation is that of a regular $m$-tuple of
weights, and after adding the extra ramification at $q_1$ to the
first component (at this step we have seen, via Theorem 4.2.1 of
[2], that automorphy propagates well) we still are in a regular
situation. We will apply a trick based on results from
 [1] to reduce to a $C$-very spread situation. \\
We choose a prime $p$ larger than all weights, larger than $t_1$ and
$q_1$, and such that the $m$ residual mod $p$ representations have
all large images, are not isomorphic to each other not even up to
twist (this, as in the previous section, can be done due to the main
theorem of [20], and regularity plus Fontaine-Laffaille theory for
the non-isomorphic condition). Recall that in
 this situation, it follows from results of [16] that the residual image of
 the tensor product is adequate, and a similar statement also holds after restriction
 to any solvable extension. We consider the mod $p$ residual representations attached
  to $f_2, f_3, ...., f_m$ and, applying Hida theory if we are in a
  residually ordinary case (which implies $p$-adic ordinariness, since we are in a
   Fontaine-Laffaille situation) or Lemma 4.1.19 in [1] combined with Theorems
    4.3.1 and 4.2.1 in [2] in case the residual representation is irreducible,
    we construct a crystalline potentially diagonalizable lift of each of these
    mod $p$ representations having arbitrarily large weight $k'_i$. We explain
    this in more detail: in the ordinary case, we consider the Hida family
    containing the $p$-adic representation attached to $f_i$ and we know that
     it contains crystalline members of arbitrarily large weight. Being ordinary
     and crystalline, any such member is potentially diagonalizable at $p$. In the
      residually irreducible case, the result from [1] is local, ensuring the
      existence of a potentially diagonalizable crystalline lift of arbitrarily
      large weight of the local at $p$ residual representation, but then
      Theorem 4.3.1 in [2] ensures the existence of a global lift with this local
      behavior\footnote{actually, to apply Theorem 4.3.1 we should replace $\Q$ by
       a quadratic imaginary field $F$ where $p$ and all primes in the level are
       split, but there is a standard argument to show that, once the lift is
       constructed over $F$, then there also exists a lift over $\Q$: the idea is
        that for the restriction of the given representation to $G_F$ one can prove,
         as a consequence of automorphy and properties of Hecke algebras, that the
         universal deformation ring $R'_S$ with the given local properties is a
         finite module over the ring of integers of certain extension of $\Q_p$,
         and by comparing this ring with the corresponding universal deformation
         ring $R_S$ for the given representation of $G_\Q$, one deduces that the same
         property holds for the latter. On the other hand, by computing the local
         deformation rings plus some cohomological calculations, one can show that
          the Krull dimension of $R_S$ is at least $1$, and this together with
          finiteness imply the existence of a $\bar{\Q}_p$ point in this ring,
           giving the lift we wanted. All this procedure holds in certain generality,
            including the case where the local conditions correspond to fixing irreducible
            components of the local universal deformation rings, thus covering our case.
For references, see [8], section 2, [17], section 4, [5], lecture 5
and [14], section 7. In a primitive form,
 the trick of comparing $R'_S$ to $R_S$ originated during the proof of Serre's conjecture, and, previous to that,
  the lower bound for the Krull dimension of the ring in works of B\"ockle}, and moreover,
  and this is also important, it can be taken of the same type as the $p$-adic representation
   attached to $f_i$ at all other places. Finally, the A.L.T. 4.2.1 of loc. cit. allows us to
    deduce modularity of the lift just constructed from the modularity of the given representation,
     the one attached to $f_i$, since the latter is known to be Fontaine-Laffaille, thus potentially diagonalizable at $p$.  \\
Let us record what we have just done in a Lemma, which we will state in sufficient generality for our
 future uses of it. For some of the Theorems we have mentioned in the above paragraph, some
 conditions on the size of the residual image are required, and for simplicity (maybe this
  condition can be relaxed) we also assume that $p > 5$ (as required in Theorem 4.3.1 of [2]) :

\begin{lema}
\label{teo:tobystrick} Let $f$ be a newform of weight $k$ and let $p
> 5$ be a prime such that the $p$-adic Galois representation
attached to $f$ is crystalline and Fontaine-Laffaille (in
particular, it can be Barsotti-Tate) locally at p. In other words,
$p > k$ and $p$ is not in the level of $f$. Suppose that the image
of the mod $p$ Galois representation attached to $f$ is large. Let
us call $\bar{r}$ this residual representation. Then, for any given
integer $E$, there is an integer $k' > E$ such that $\bar{r}$ has a
$p$-adic lift which is modular, associated to a modular form $f'$ of
weight $k'$, and this $p$-adic lift is crystalline and potentially
diagonalizable at $p$ (in particular, $p$ is not in the level of
$f'$). Moreover, we can choose $f'$ such that the $p$-adic Galois
representations attached to $f$ and $f'$ connect locally at any
prime $\ell \neq p$ (they have the same types).
\end{lema}

Since we can take them to be arbitrarily large, we declare that the ``new" weights $k'_2, ....k'_m$, which are chosen one at a time, are going to satisfy the inequalities in the definition of $C$-very spread. Thus we have reduced to a situation where the weights are $C$-very spread (we recall that we have fixed $C = q_1 + 1$).\\
We have to check that, in the above process, while moving all components from the second to the last through suitable congruences, some A.L.T. can be applied to the $m$-fold tensor products guaranteeing that automorphy of one side is equivalent to automorphy of the other. But we have Fontaine-Laffaille (thus potentially  diagonalizable) components on one side, and we have a $C$-very spread set of potentially  diagonalizable components on the other side and we also have already remarked that the adequacy condition on the residual image of the tensor product is satisfied, thus Theorem 4.2.1 in loc. cit. (or, if we don't want to assume that $p$ is larger than $2 (2^m + 1)$, we can use the improvement in [11]) can be applied and we are fine. \\

To ease notation, let us call these new weights again $k_2, ...., k_m$, in the understanding that we now know that, together with $k_1$, this is a $C$-very spread $m$-tuple of weights.

\section{In which the other Good-Dihedral primes enter the picture}
Before going on, and in order to guarantee that the components from the second to the last will
 also have large residual images in what follows, we are going to add Good-Dihedral primes to all
  of them. This is done exactly as we did for the first component in section 3, with the following conditions: \\
We apply this procedure one component at a time. Suppose we have
applied it to the components from the second to the $j-1$-th one,
then for the $j$-th component we choose a bound $B_j > C = q_1 + 1$,
and we pick a prime $t_j > B_j$ such that the required conditions on
all residual mod $t_j$ images
 are satisfied (see section 3 for the case of the first component and do exactly the same), larger
  than the $m$ weights $k_i$ and different from any of the $t_i, q_i$ $i= 1,..., j-1$ already introduced,
  and modulo this prime $t_j$ introduce a Good-Dihedral prime $q_j$ such that it is different to all the
  $q_i$ and $t_i$, $i=1,2,..., j-1$, and such that it satisfies the definition of Good-Dihedral prime with respect to the bound $B_j$. \\
Since all the bounds $B_j$, $j >1$, are larger than $C$, as long
 as we work in characteristic $p \leq C$, as we would do while
  moving through the skeleton chain in the first component, we
   know that the residual images of all but the first component
    will be large (and the residual image of the first component is also large\footnote{with the exception of the very last move, which will be discussed in Section 9}
    because of how the skeleton chain was built, see section 3 and [10]). This is
    even true if the characteristic was equal to $3$ or $5$ (see the Introduction of [10] for the
    definition of large in this case) because we have an element of large order $t_i$ in
    the projective residual image. \\
Observe also that the residual representations from now on will not be isomorphic to each
other, not even up to twist, because of their differences in ramification, this is why we
insisted to pick all the $q_i$ different, and this remains true until the end since the
 Good-Dihedral primes $q_2, ...., q_m$ will never be removed. \\
This being said, and recalling once again that largeness is preserved by restriction to solvable extensions, we see (this follows from [16], see footnote in section 3) that as long as we work in characteristics $p \leq C$ the condition of adequacy of the residual image of the $m$-fold tensor product after restriction to the cyclotomic extension, required for A.L.T. to work, will hold. \\
We will work in characteristic $p \leq C$ most of the time, but
 at a few steps we will also need to work in larger characteristics:
  we will choose them different from the $t_i , q_i $, so at least the
   local behavior at $q_i$ will not be altered. When working in such a
    large characteristic $p$, we will have some freedom to choose it,
    thus largeness of the residual images will again be easily established,
    usually through Ribet's large image theorem. \\

\section{In which we show that at least in most of the steps things work well}
We have reduced the proof of the main theorem to a situation in which the $m$ components have each
 its own Good-Dihedral prime $q_i$, the weights are $C$-very spread, and now we want to keep moving
 the first component through the skeleton chain (the first step, which is the addition of a
 Good-Dihedral prime, has already been done). We are not going to explain this one move at a time,
  we are just going to classify the three possible type of congruences showing up in the
  skeleton chain, and we will see how to proceed in each case. We begin in this section,
   with the most common type of congruence (this one occurs in all but a few moves inside
    the skeleton chain, where a few means approximately $9$). Let us define it:\\

{\bf Definition}: We say that a mod $p$ congruence between two cuspforms is of type (A) if both of the forms are potentially diagonalizable at $p$, of different weights, and they are of the same type at any other prime. \\

Typically, type (A) in the skeleton chain appears in congruences between Fontaine-Laffaille and potentially Barsotti-Tate representations. But there are other cases, like in one of the final four moves (see section 3). \\

Let us see that if the congruence in the first component is of type (A), while the other $m-1$ components are left unchanged, then Theorem \ref{teo:enlamenor} can be applied to show that automorphy propagates well through this congruence of $m$-fold tensor products. Here we are taking the $2$-dimensional Galois representations in the first component as $R$ and $R'$, and $S= S'$ to be the $m-1$-fold tensor product of the other components. Since we have that the two weights $k_1$ and $k'_1$ involved in the congruence in the first component are both bounded by $C$, and moreover that the weights are $C$-very spread, we see that the extra-regularity condition (4) in Theorem \ref{teo:enlamenor} is satisfied, as explained after the definition of $C$-very spread (see section 4). \\
As for the other three conditions in Theorem \ref{teo:enlamenor}, condition (1) is satisfied since we are in a type (A) case, condition (2) is automatic since $S = S'$, and for condition (3) this follows from known properties of the skeleton chain in the first component plus the Good-Dihedral primes $q_i$, $i >1$, in the other components, as explained at the end of section 5. \\
Thus, we conclude that automorphy propagates well if the congruence in the first component is of type (A).\\

Remark: There is one and only one congruence in the skeleton chain with residual characteristic smaller than $7$, and it is the mod $3$ congruence in one of the last four moves introduced in Section 3. As we explained there, this congruence is between potentially diagonalizable representations of different weights, and the type at primes other than $3$ is not altered, so it is of type (A). The residual image is large because it contains $\SL_2(\F_{27})$, and for the components from the second to the last as we have already remarked the existence of Good-Dihedral primes also forces the residual images to be large (recall that large in characteristic $3$ means that it must contain $\SL_2(k)$ with $k$ a finite extension of $\F_3$ different from $\F_3$). Thus Theorem \ref{teo:enlamenor} applies in this congruence. \\

\section{What do to in cases where no A.L.T. would apply: the ``genetic manipulation" of the other components}

There are other two kind of congruences showing up while moving along the skeleton chain, we define them now: \\

{\bf Definition}: We say that a mod $p$ congruence between two cuspforms is of type (B) if both of the forms are ordinary at $p$, of different weights, and they are of the same type at any other prime. \\

{\bf Definition}: We say that a mod $p$ congruence between two cuspforms is of type (C) if both of the forms are Barsotti-Tate at $p$, and they are of different type at some prime $\ell \neq p$. \\

Remarks: \\
--1) As the attentive reader may have observed, the congruence in Step 1, as done in section 3, to add a Good-Dihedral prime, is not of any of the three types defined. We simply exclude it from our considerations, since we have already seen that an A.L.T. was applicable there (alternatively, we could have given a slightly more general definition of type (C)).\\
--2) In type (B), we can assume that at least one of the two ordinary representations is not potentially crystalline, else this will be type (A). \\
--3) Whenever a congruence of type (B) or (C) occurs in the skeleton chain, we have $p \geq 7$ (see the remark at the end of section 6).\\

The problem is that if the congruence in the first component is of any of these two types, and the other components are not known to be potentially diagonalizable (the typical situation is that the prime $p$ is smaller than some of the weights, and the representations are crystalline but not Fontaine-Laffaille), no A.L.T. can be applied to propagate automorphy. So, what can we do to keep moving forward in our proof? \\

What we propose is to change, through suitable congruences (where suitable A.L.T. will guarantee that automorphy propagates well) all components from the second to the last one, one at a time, in order to make them potentially diagonalizable at $p$. This is what we informally call a ``genetic manipulation". AFTER FINISHING WITH THIS PROCEDURE, we will be able to go through the move in the skeleton chain, i.e., in the first component, either involving type (B) or type (C), as we will explain later. \\

First of all let us stress the key idea that has allowed us to do this: assuming that we are in the $i$-th component,
 $i > 1$, and the $p$-adic Galois representation there is in a Fontaine-Laffaille situation (we will show how to
 reduce to such a case, ONE COMPONENT AT A TIME:  the pace is key, else this would be impossible if $p$ is smaller,
  say, than $m$, because of regularity), we can consider, using Lemma \ref{teo:tobystrick}, a congruence with  another modular form of arbitrarily large weight $k'_i$ such that it is crystalline and potentially diagonalizable at $p$. \\

The procedure is a bit painful, since we need to make sure that suitable A.L.T. do apply at each move performed during the manipulation (and this implies keeping track of local conditions, plus regularity, plus residual images), and we also want to make sure that when we finish we haven't altered a fundamental fact: we start with an $m$-tuple of weights $k_1 , k_2 , ....., k_m$ which is $C$-very spread, and the output must be another $m$-tuple of weights (we do not change the first one) $k_1, k'_2, k'_3, ......., k'_m$ also $C$-very spread. In this procedure, we will also need to introduce auxiliary primes (larger than the bounds $B_i$) in the level of all the components from the second to the last one, and we also have to make sure that this does not alter the fact that their Good-Dihedral primes $q_i$ remain ``Good-Dihedral primes with respect to the bound $B_i$'' after this extra ramification has been introduced. \\

So, here is where the hard work begins. We will divide the algorithm in two cases, not depending on whether
 we are acting ``just before" a congruence of type (B) or of type (C) in the skeleton chain,
  but depending on the following fact: Let $k_1$ be the input weight in the skeleton
   chain, i.e., we have a newform of weight $k_1$ and we want to move through a congruence with another newform, possibly of different output  weight. Then we will consider separately the case $k_1 =2$ and the case $2 < k_1 \leq C$. \\

Remark: To avoid any confusion with the notions of input and output
weights, we insist
 that we are moving along the skeleton chain starting on some level 1 newform $f_1$ and finishing in a small level CM form. Thus, if we were to draw this chain from left to right, the leftmost point will correspond to $f_1$, and the input of a congruence is just at the left of its output. In the end, automorphy propagates just the other way round, from right to left, this is why I wanted to clarify these somehow arbitrary notions.\\

We will assume first that the input weight of the first component is
not $2$, the easiest case, even if this case rarely occurs. Let us
see how to transform, one at a time, all other components into
representations that not only still have very large weights but also such that they are potentially diagonalizable at $p$. This is done simply, for example to
the second component,  by first taking a weight $2$ lift modulo any
sufficiently large prime $p_2$. This prime must be taken to be, and
this we have seen already a couple of times how can be done with
standard techniques, such that the residual image of the $m$-fold
tensor product satisfies the usual adequacy condition  (after
restriction to the cyclotomic field) and larger than all weights, in
order to be in a Fontaine-Laffaille situation. We also take it to be
different from any other prime ever considered thus far in this
paper. Observe also that after this step the $m$-tuple of weights
become (in increasing order):
$$ (2, k_1, k_3 , ....., k_m)$$
which is clearly spread since $k_1 > 2$. Automorphy propagates well through this congruence, due to Theorem 4.2.1 of [2], since both potentially Barsotti-Tate and Fontaine-Laffaille cases are potentially diagonalizable. \\
Just one detail before moving on: we must choose the prime $p_2$ to be a square mod $q_2$,
 so that $q_2$ is a square mod $p_2$. This way, even if we have introduced extra ramification
  at $p_2$ in the second component, the Good-Dihedral prime $q_2$ will still be enough to ensure
   large residual images for the second component, in characteristics up to $B_2$, as can be easily checked (using the usual method for proving that Good-Dihedral primes give large residual images). \\
Now we are ready for the key step: we have already reduced the second component to a weight $2$ situation, so that if we move to characteristic $p$, $p$ being the prime where the congruence in the skeleton chain will take place, we are in a Barsotti-Tate situation ($p$ was not, and is not, in the level of the newform in the second component). So now we apply Lemma \ref{teo:tobystrick} to the mod $p$ representation in the second component, since we know that $p>5$ (see Remark 3 above), and the residual image is large (due to the Good-Dihedral prime $q_2$ and $p < C < B_2$) it can be applied and it gives us another modular lift, of weight $k'_2$ arbitrarily large, crystalline and potentially diagonalizable at $p$, and preserving the types locally at all other ramified primes.\\
We pick $k'_2 > 2 C + k_3 + ..... + k_m $. In this mod $p$ congruence in the second component, we have as input $m$-tuple of weights: \\
$$(2, k_1, k_3, ........, k_m)$$
And the output weights are:
$$(k_1, k_3, ............, k_m, k'_2)$$
Then, the above inequality is enough to guarantee that the $m+1$-tuple:
$$ (2, k_1, k_3, ......, k_m, k'_2) $$
is spread. Here we are using the fact that the initial weights $(k_1, k_2, ......, k_m)$ were $C$-very spread
(and in particular $k_1 < C$) and the inequality imposed on $k'_2$. \\
Thus we can apply Theorem \ref{teo:enlamenor} to see that modularity propagates well through
this congruence of $m$-fold tensor products. We are taking $S = S'$ to be the $m-1$-fold tensor
 products of all but the second component, and $R$ and $R'$ the modular Galois representations in the congruence
  of the second component. We have just checked that the extra-regularity condition (4) is satisfied, conditions (1) to (3) are easy to check (as explained in section 5, since we are in characteristic $p < C$, the adequacy condition on the residual image holds). \\
The rest is just an iteration of the same idea: after modifying the
$i-1$-th component, for any $2 < i \leq m$, we manipulate the $i$-th
component in a similar way: we choose as above a sufficiently large
prime $p_i$, larger than all weights (in increasing order)
$$(k_1, k_i, ..... k_m, k'_2, ....k'_{i-1} )$$ and satisfying the
same conditions we asked to $p_2$, but now translated to the $i$-th
component, in particular we ask this prime to be a square mod $q_i$.
We take a weight $2$ lift in the $i$-th component modulo $p_i$. Here
we check again that Theorem 4.2.1 of [2] allows automorphy to
propagate. Then, we move to characteristic $p$, and again via Lemma
\ref{teo:tobystrick} we make a congruence with a modular form of
weight $k'_i$, which we can and will choose to satisfy:
$$ k'_i > 2 C + k_{i+1} + ....... + k_m + k'_2 + ...... + k'_{i-1} \qquad \qquad (@)$$
With this it is easy to recursively show that in this last congruence, the extra-regularity condition (4) in Theorem \ref{teo:enlamenor} holds, in other words, that the $m+1$-tuple:
$$ (2, k_1, k_{i+1}, ......, k_m, k'_2, ......, k'_i ) $$
is spread. Conditions (1) to (3) of this A.L.T. are again easy to check, thus we conclude that automorphy propagates well. \\
After applying the algorithm to the $m$-th component, we end up with an $m$-tuple of newforms of weights:
$$(k_1 , k'_2, .......k'_m )$$
such that all components but the first one are crystalline and potentially diagonalizable at $p$. Moreover, by the way the weights $k'_i$ were chosen it is easy to see that this new set of weights is again $C$-very spread. \\
Thus we can proceed and go through the mod $p$ congruence in the first component, either if it is of type (B) or of type (C). \\
We can apply either Theorem \ref{teo:mixed} (for type (B)) or Theorem 4.2.1 in [2] and its variant in [11] (for type (C)), because we will either be in a ``mixed" situation such as the one Theorem \ref{teo:mixed} was meant to dealt with, or in a potentially diagonalizable situation.  The weights are $C$-very spread, so no matter the output value of the weight in the first component (it is always bounded by $C$) we know that we have a congruence between two regular $m$-fold tensor products, and since $p < C$ we also know that the residual image is sufficiently large.\\

Now we have to deal with the other case, the case with $k_1 = 2$, which is a bit more complicated.
 All auxiliary primes considered will be assumed to be different between them and different from
 any prime previously considered in this paper.\\
We again have to alter $k_2, k_3, ...., k_m$ and we will do it one
at a time, so to begin with let us ``transform" the second component
into a representation of some large weight $k'_2$ that
 is potentially diagonalizable at $p$.
The moves are as follows (recall that $p$ is at least $7$): starting with weights $(2, k_2, ...., k_m)$,
 modulo some prime $u_2$ larger than all weights such that the first component is ordinary,
  and residual images are all large we change the weights to $(3, k_2, ....., k_m)$, using Hida
   theory and specializing to weight $3$. Observe that such a prime exists because of Ribet's
   large images theorem and the fact that weight $2$ newforms have a density one set of ordinary primes (see [13]).
   Observe also that A.L.T. 4.2.1 in [2] applies to this congruence of $m$-fold tensor products
    where we are only changing the first component,  because we know that the residual image is
     sufficiently large, the first components are potentially diagonalizable (the weight $3$ modular form is
      potentially crystalline and ordinary at $u_2$), the other components are Fontaine-Laffaille,
       thus also potentially diagonalizable, and of course $3 < C$ so the new set of weights is also $C$-very spread. \\
It is perhaps worth remarking at this point that we have violated a sacred rule: we have introduced
extra ramification, at $u_2$, in the first component, thus altering the skeleton chain. To remedy this,
 we are going to undo this move as soon as possible, just in two more moves, so to rejoin the skeleton
 chain after this small detour, safe and sound. \\
The next move is modulo some large prime $v_2$,  here we take  a weight $2$ lift in the second component,
 i.e., we change weights to  $(3, 2, k_3, ..... k_m)$. This set is easily seen to be spread. This move
 is as the one we did in some large characteristic while treating the previous case, and again we take
  this prime to be larger than all weights, so that all residual images are large, and such
  that $v_2$ is a square mod $q_2$. We easily check that A.L.T. 4.2.1 in [2] applies to
  this congruence. \\
Then, modulo some prime $w_2$ larger than all weights, such that the second component is
 ordinary, all residual images are large, and such that $w_2$ is a square mod $q_2$
 (we use
 again the fact that the set of ordinary primes has density $1$ for a weight $2$ newform to see that such a prime exists), using Hida theory we change weights
 to $(3, 6, k_3, .... , k_m)$, which are spread. As with the prime $u_2$, we check
 that every two-dimensional representation involved in this congruence is potentially diagonalizable at $w_2$,
  and we easily see that A.L.T. 4.2.1 in [2] applies.
Then, we move back to characteristic $u_2$,  and we undo the move on the first component:
 we obtain weights $(2, 6, k_3, ......., k_m)$, which are spread. Observe that the residual
  image of the first component is large for the very simple reason that it was so two moves
   ago (this was one of the conditions to choose this prime then) and we haven't changed the
    first component since then. In fact, the same applies to the residual images of all but
     the second component: they were large two moves ago, and the newforms haven't changed,
      so it is clear that their mod $u_2$ images are still large. But the second component
      has changed, and thus a priori $u_2$ could be in the exceptional set for Ribet's theorem
       for the weight $6$ newform that is now in the second component. But we
       can check using the techniques in [12] that this is not the
       case, that the residual image in the second component is also
       large. What we do is to combine the Good-Dihedral prime with a
       local analysis at $u_2$: because of the Good-Dihedral prime
       $q_2$, applying its existence in the usual way, we see that the residual image can fail to be large only if it is
       dihedral and such that the quadratic field where it becomes reducible
       is ramified at $u_2$.
        But looking at the action of the tame
       inertia group at $u_2$, since we are in a Fontaine-Laffaille
       situation, we can apply ideas of Ribet (which are recorded
       for example as Lemma 3.1 of [10]) to deduce that such a
       dihedral case can only occur if $u_2 = 2k -1$ or $u_2 = 2k -
       3$, where $k$ is the weight of the given newform. Since in
       our case $k=6$ and $u_2 > 12$, we conclude that the residual
       image is large. \\
             Having checked what was required on residual images and since all two-dimensional representations involved
              are potentially diagonalizable, we see that A.L.T. 4.2.1 in [2] applies to this congruence. \\
Now we can move to characteristic $p$ and apply Lemma \ref{teo:tobystrick}  to the second
 component: since $p \geq 7$, the $p$-adic representation in the second component is
 Fontaine-Laffaille, so we can make a congruence with a potentially diagonalizable crystalline modular
  representation of arbitrarily large weight $k'_2$. Thus the output weights are
  (in increasing order) $(2, k_3, ....., k_m, k'_2)$. It is important to stress
   that the newform of weight $2$ in the first component is exactly the one we started with at
    the beginning of the algorithm, because the only change we have performed to the first
    component was soon undone. Thus, it belongs to the skeleton chain, and we know
    (this is part of the structure of the skeleton chain) that since the next move
     in the skeleton chain is in characteristic $p$, the mod $p$ residual image
     of this representation is large (there is an exception to this, at the last
      step of the skeleton chain where the residual image is dihedral, but we are
       going to dedicate Section 9 to explain what are the differences in that case).
       As for the residual images of the other components, we know since $p < C$ that they are large.
        Let us specify the lower bound for $k'_2$: we take it exactly as in the other case, i.e., we
         declare that $k'_2 > 2 C + k_3 + ..... + k_m $. With this, it is easy to see that the extra-regularity
          condition (4) in Theorem \ref{teo:enlamenor} holds, namely, that the $m+1$-tuple of weights:
$$ (2, 6, k_3 , ......., k_m, k'_2) $$
is spread (note the similarity with the previous case). \\
Conditions (1) to (3) of that A.L.T. are easily seen to hold, since we have already
 explained that residual images of all components are known to be large (and due to their
 differences in ramification, these residual representations are not isomorphic to each other,
 not even up to twist) thus we conclude that it applies and ensures that automorphy propagates
 well through this congruence. \\
The rest is just an iteration of the above procedure and it works well for the same standard
 arguments. To make the moves in the $i$-th component, for any $i > 2$, assuming we have
 already done so for the previous components, one repeats similar moves with large primes
  $u_i$, $v_i$, $w_i$, the last two of them congruent to squares mod $q_i$ and satisfying
  properties similar to those required to $u_2, v_2, w_2$, then modulo $u_i$ the first move
  is undone, and one culminates with an application of Lemma \ref{teo:tobystrick} in
  the $i$-th component to change from a weight $6$ $p$-adic modular representation to
   a potentially diagonalizable and crystalline one of weight $k'_i$ chosen as in the previous case, see
   formula (@). \\
As happened in the second component, it is easy to check that in all these moves some A.L.T.
 applies, thus ensuring that automorphy propagates well. In particular, the extra-regularity
 condition (4) in Theorem \ref{teo:enlamenor}, a theorem that we apply in the last move in
  the $i$-th component, holds as in the previous case as follows from formula (@). \\
After applying these moves to the last component, we end up with a
set of  weights:
 $(2, k'_2, ....., k'_m)$ which is $C$-very spread, as follows from formula (@) as
 in the previous case. \\
Also, all components from the second to the last are potentially diagonalizable at $p$, and residual
images modulo $p$ of all components are large. Thus, we can now  go through the mod $p$
congruence in the first component, either if it is of type (B) or of type (C), by applying
a suitable A.L.T., either Theorem \ref{teo:mixed} or the A.L.T. in [11]. \\

We have thus been able, through some manipulations, to reduce to a situation where even if
the congruence in the skeleton chain is of type (B) or (C), a suitable A.L.T. can be applied.
The modifications we have performed in the components from the second to the last are harmless:
 on one hand, weights have changed, but they still are $C$-very spread, thus regularity in the
 following steps is ensured; on the other hand,  some extra ramification has been introduced
  into them, but in such a way that the corresponding Good-Dihedral prime $q_i$ ($i > 1$)
  still plays the same role, thus in particular ensuring that when working in any characteristic
   $p < C$ the corresponding residual image will be large. \\
Since we have explained in the previous section that in the type (A) case a suitable A.L.T.
applies, this means that with these manipulations we have been able to go through all congruences
 in the skeleton chain in such a way that a suitable A.L.T. applies at each move, thus reducing
  the proof of the main theorem (which is Theorem \ref{teo:main}) to a case where one component
   is a CM form. But since automorphy is preserved by tensoring by a CM form (this is one of
   the key observations that makes Harris' trick work) this means that we have reduced the proof
    to a case of a regular tensor product of $m-1$ modular Galois representations. In  section
    10 we will explain that despite the fact that the newforms that we found in this $m-1$-fold
    tensor product are no longer of level $1$, still we can conclude the proof of the main theorem
    by induction on $m$. Before that, we dedicate two sections to clarify a couple of technical points,
     for the reader's convenience. But it is fair to say that we have already done most of the job,
     what the three following sections contain are just minor issues.

\section{Remarks on Galois conjugation, and on twists}
In the skeleton chain, there are two minor points that should be clarified: \\
Twists: sometimes congruences are up to twist by a finite order character, i.e., we are changing the output representation by a twist of it in order, for example, to minimize the residual Serre's weight. When we put both representations inside  $m$-fold tensor products, the fact that we are working up to twist is harmless because automorphy is preserved by twists and all conditions in all the A.L.T. are also preserved by twisting (in particular, potential diagonalizability). \\
Galois conjugation: finally, in the skeleton chain we also sometimes replace a modular form by a Galois conjugate of it, in particular, this is why the fact that the ``bottom space" contains a single orbit is good enough to use it as a bridge to propagate automorphy. But when we put the modular form (as first component) inside an $m$-fold tensor product, automorphy will in principle not necessarily be preserved by Galois conjugation of the first component. To remedy this problem, we do the following: whenever in the first component one applies Galois conjugation, at the same time one conjugates all other components by the same Galois element (i.e., by a compatible choice of Galois elements). This way, we are inducing a conjugation on the $m$-fold tensor product, which is a valid move since it preserves automorphy. The alteration infringed in the components from the second to the last is harmless: after conjugating, they preserve their weights and their levels and their Good-Dihedral primes $q_i$, so this is an optimal way to deal with the issue of Galois conjugation.\\

\section{In which we explain why when one component has a dihedral but not bad-dihedral residual image A.L.T. still work well}

At several points, we have used the argument that the residual image of
 the $m$-fold tensor product satisfies the adequacy condition required
 to apply certain A.L.T. because the residual images of all components
 were large and not isomorphic to each other, not even up to twist. \\
This is true at all congruences between $m$-fold tensor products considered
 in this paper, except for the case where we get to the end of the skeleton
  chain and there we have a mod $13$ congruence between two representations,
   in the first component, with dihedral residual image. Moreover, since the
    congruence is of type (C), we have to go through the ``genetic manipulation"
    here, and this involves working in characteristic $13$ several times during
    the process of changing the weights of the components from the second to
    the last one, and this is done with a first component having dihedral residual image. \\
Thus, we need to make sure that at all these moves in characteristic $13$
the condition on the size of the residual image required to apply an A.L.T.
is satisfied, despite the fact that residual images of the components are not all large. \\
For the components from the second to the last, we know that due to their
 Good-Dihedral primes the mod $13$ residual images are large, and we can argue
 in the usual way to see that for the $m-1$-fold tensor product having them as
  components the residual image is adequate, even after restriction to the
   cyclotomic field of $13$-th roots of unity. Then, we can apply
   Lemma A.3.1 in [1] to check that tensoring with the first
   component gives a residual image that also satisfies this
    property. The only thing that has to be checked is that
    the dihedral representation that we have residually in the first
    component is irreducible when restricted to the compositum $M$ of the
     cyclotomic field of $13$-th roots of unity and the fixed field of the
      kernel of the $m-1$-fold tensor product of the other residual components.
       But this follows from the combination of two facts: first, this dihedral
        representation is induced from a character of $\Q(\sqrt{-3})$ (recall
        that the congruence modulo $13$ in the first component is with a CM
        elliptic curve of conductor $27$) and since the size of the projectivized
        residual image is larger than $4$ (this can be checked easily by looking
        at the inertia group at $13$) it will become reducible after restriction
         to a Galois number field only if such field contains $\Q(\sqrt{-3})$,
         and second, the only component that ramifies at $3$ is the first one,
         thus $M$ does not contain $\Q(\sqrt{-3})$. \\
Thus we conclude that even if in this case the residual image of the first component
 is not large, the adequacy condition on the residual image of the $m$-fold tensor product is satisfied.\\

\section{Completing the induction process: rewinding}

At the end of our chain of congruences, we have reduced the proof of automorphy of the given $m$-fold tensor
 product of level $1$ cuspform to that of an $m-1$-fold regular tensor product, but this time the components
  are newforms of level greater than $1$. In order to have a proof by induction on $m$, we should now
  ``safely connect" this tensor product with another one whose components are level $1$ newforms. \\
But what may seem complicated is rather simple: we are looking for a chain of congruences propagating
 automorphy ending in an $m-1$-fold tensor product of level $1$ newform, but we have already constructed
  such a chain! Simply consider all the moves we have done with the $m$-fold tensor products, starting
  in a ``level 1 case", through our proof, delete from all this the first component and forget about the
   congruences involving exclusively the first component, and reverse the resulting chain of congruences
    of $m-1$-fold tensor products. This is clearly a chain of $m-1$-fold tensor products starting where
     we wanted to start and ending, as we wanted, in a level $1$ case. Is it safe, i.e., do A.L.T. apply
      through this chain (as usual, we need to apply them the other way round, i.e., starting in the level $1$ tensor product)? The answer is yes, and this is due to the following facts: \\
1) The chain of $m$-fold tensor products that we have built in previous sections is safe. \\
2) Moreover, at each step we have applied A.L.T. that work in both directions, so that they can be used to propagate automorphy not only from right to left, but also from left to right (this allows us to reverse now). \\
3) For all the A.L.T. that we have applied to $m$-fold tensor products, we have checked the required local conditions on each component (because these are conditions that behave well under tensor products) and we have also checked that the residual image of the tensor product was sufficiently large by checking that the images of the components were so (and that there were no unwanted isomorphisms between them). As for regularity, it is clear that if it holds for the case of $m$ components it will also hold after dropping the first component. The same applies to the extra-regularity condition (4) in Theorem \ref{teo:enlamenor}. Thus we conclude that by dropping the first component we are not altering the fact that the chain is safe. \\

This concludes the reduction of the proof to the case of an $m-1$-fold regular tensor product of level $1$ newforms. Thus, since the claim is trivial for $m=1$,  the proof of the main theorem follows by induction on $m$.

\section{Bibliography}

[1] Barnet-Lamb, T., Gee, T., Geraghty, D., {\it    Serre weights
for rank two unitary groups}, preprint; available at www.arxiv.org
\newline
[2]   Barnet-Lamb, T., Gee, T., Geraghty, D., Taylor, R., {\it
Potential automorphy and change of weight}, preprint; available at
www.arxiv.org
\newline
[3] Berger, L., Li, H.,  Zhu, H., {\it Construction of some families
of $2$-dimensional crystalline representations},  Math. Annalen {\bf
329} (2004), 365-377
\newline
[4] Billerey, N., Dieulefait, L., {\it Explicit Large Image Theorems
for Modular Forms}, preprint; available at www.arxiv.org
\newline
[5] B\"ockle, G., {\it Deformations of Galois representations}, in
``Elliptic curves, Hilbert modular forms and Galois deformations",
H. Darmon, F. Diamond, L. Dieulefait, B. Edixhoven, V. Rotger
(eds.),
 Progress in Math., Birkhauser, to appear.
\newline
[6] Calegari, F., {\it   Even Galois Representations and the
Fontaine-Mazur Conjecture. II }, JAMS {\bf 25} (2012), 533-554
\newline
[7] Caraiani, A., {\it  Monodromy and local-global compatibility for
$l=p$}, preprint, available at www.arxiv.org
\newline
[8] Clozel, L., Harris, M., Taylor, R., {\it   Automorphy for some
$\ell$-adic lifts of automorphic mod $\ell$ representations}, Pub.
Math. IHES {\bf 108} (2008), 1-181.
\newline
[9] Dieulefait, L., {\it Langlands base change for $\GL(2)$}, Annals
of Math. {\bf 176} (2012) 1015-1038
\newline
[10] Dieulefait, L., {\it Automorphy of $\Symm^5(\GL(2))$ and base
change}, preprint; available at www.arxiv. org
\newline
[11] Dieulefait, L., Gee, T., {\it Automorphy lifting for small
$l$}, Appendix B to [10]
\newline
[12] Dieulefait, L., Wiese, G., {\it    On Modular Forms and the
Inverse Galois Problem}, Trans. AMS {\bf 363} (2011),  4569-4584
\newline
[13] Gee, T., {\it The Sato-Tate conjecture for modular forms of
weight 3}, Doc. Math. {\bf 14}  (2009), 771-800
\newline
[14] Gee, T., Geraghty, D., {\it   Companion forms for unitary and
symplectic groups},  Duke Math. J. {\bf 161} (2012),  247-303
\newline
[15] Guralnick, R., {\it  Adequacy of representations of finite
groups of Lie type},  Appendix A to [10]
\newline
[16] Guralnick, R., Herzig, F., Taylor, R., Thorne, J., {\it
Adequate subgroups}, Appendix to  ``On the automorphy of $l$-adic
Galois representations with small residual image" by J.Thorne, J.
Inst. Math. Jussieu {\bf 11} (2012),    855-920
\newline
[17] Kisin, M., {\it     Modularity of 2-dimensional Galois
representations}, Current Developments in Mathematics 2005, 191-230
\newline
[18] Khare, C., Wintenberger, J-P., {\it Serre's modularity
conjecture (I)}, Invent. Math.  {\bf 178} (2009) 485-504
\newline
[19] Ramakrishnan, D., {\it Modularity of the Rankin-Selberg
$L$-series, and multiplicity one for $\SL(2)$}, Annals of Math.
{\bf 152} 45-111 (2000)
\newline
[20] Ribet, K., {\it On $\ell$-adic representations attached to
modular forms. II}, Glasgow Math. J. {\bf 27} (1985) 185-194

\end{document}